\documentclass[12pt]{amsart}
\usepackage{amssymb}
\usepackage{graphicx}
\usepackage{amsmath}
\usepackage{amsfonts}

\newcommand{\rsphere}{ {\mathbb C_\infty} }
\newcommand{\ucirc}{ {\mathbb S^1}}

\newcommand{\bd}{\partial}
\newcommand{\Int}{\operatorname{Int}}

\newtheorem{thm}{Theorem}
\newtheorem{lem}[thm]{Lemma}
\newtheorem{cor}[thm]{Corollary}
\newtheorem*{thma}{Theorem}
\theoremstyle{definition}
\newtheorem{ques}[thm]{Question}
\newtheorem{defn}[thm]{Definition}
\newtheorem*{mconj}{Makienko's Conjecture}

\theoremstyle{remark}

\begin{document}

\title[Makienko's Conjecture and Indecomposability]{Any counterexample to Makienko's conjecture is an indecomposable continuum}

\author[C.~P.~Curry]{Clinton P.~Curry}
\email[Clinton P.~Curry]{clintonc@uab.edu}

\author[J.~C.~Mayer]{John C.~Mayer}
\email[John C.~Mayer]{mayer@math.uab.edu}

\address[Clinton P.~Curry and John C.~Mayer]
{Department of Mathematics\\ University of Alabama at Birmingham\\
  Birmingham, AL 35294-1170}

\author[J.~Meddaugh]{Jonathan Meddaugh}
\email[Jonathan Meddaugh]{jmeddaugh@math.tulane.edu}

\author[J.~T.~Rogers,~Jr.]{James T. Rogers, Jr.}
\email[James T. Rogers, Jr.]{jim@math.tulane.edu}

\address[Jonathan Meddaugh and James T.~Rogers,~Jr.]
{Department of Mathematics\\ Tulane University\\
  New Orleans, LA 70118}

\keywords{indecomposable continuum, Makienko's conjecture, Makienko conjecture, Julia set, holomorphic
dynamics, invariant Fatou component, complex dynamics, buried point, residual
Julia set}

\subjclass[2000]{Primary: 37F20; Secondary: 54F15}

\thanks{We thank the Departments of Mathematics at Tulane University and Nipissing University (Ontario), and the Fields Institute at the University of Toronto for the opportunity to work on this paper
in pleasant surroundings.}

\begin{abstract}
  Makienko's conjecture, a proposed addition to Sullivan's dictionary, can be stated as follows: 
  The Julia set of a rational function $R:\rsphere \rightarrow \rsphere$ has buried points if and only if no component of the Fatou set is completely invariant under the second iterate of $R$.
  We prove Makienko's conjecture for rational functions with Julia sets that are decomposable continua.
  This is a very broad collection of Julia sets; it is not known if there exists a rational functions whose Julia set is an indecomposable continuum.
\end{abstract}

\date{\today}

\maketitle

\centerline{Dedicated to Bob Devaney on the occasion of his 60th birthday.}
  
\begin{section}{Introduction}

Let $R:\rsphere\to\rsphere$ be a rational function, where $\rsphere$ denotes the Riemann sphere.
The \emph{Fatou set} of $R$, denoted $F(R)$, is the domain of normality for the family of functions $\{R^i \, \mid \, i \in \mathbb N\}$.
A component of the Fatou set is called a \emph{Fatou component}.
The \emph{Julia set} of $R$, denoted $J(R)$, is the complement of $F(R)$.
In the case that the degree of $R$ is not one, the Julia set is a non-empty, compact, perfect subset of $\rsphere$.
It is well-known that $J(R^n)=J(R)$ and $F(R^n)=F(R)$ for any integer $n \ge 1$.

A set X is said to be \emph{completely invariant} under $R$ provided that $R(X)=X=R^{-1}(X)$.
Both $J(R)$ and $F(R)$ are completely invariant under $R$. 
In fact, $J(R)$ is the smallest closed subset of $\rsphere$ which is completely invariant under $R$ and contains at least three points \cite{Beardon:1991}.
On the other hand, the Fatou set may have large subsets that are completely invariant under $R$.
For example, if $R$ is a polynomial, then the basin of attraction of infinity is a completely invariant component of the Fatou set.
If $R(z)=z^{-2}$, then the Fatou set of $R$ has no completely invariant component, while the Fatou set of $R^2(z) = z^{4}$ has two completely invariant components.

A point of the Julia set is said to be \emph{buried} if it does not belong to the boundary of a Fatou component.
The set of all buried points of the Julia set of a map is called the \emph{residual Julia set}.
McMullen \cite{McMullen:1988} presented the first rational maps with non-empty residual Julia sets.
He showed that functions of the form $z \mapsto z^n+\lambda/z^d$ have Julia sets homeomorphic to the product of a Cantor set and a Jordan curve whenever $1/n+1/d<1$ and $\lambda$ is sufficiently small.
In this case, there are uncountably components of the Julia set which do not intersect the boundary of any Fatou component.
Later, John Milnor and Tan Lei \cite{Milnor:1993} and Devaney, Look, and Uminsky \cite{Devaney:2005} exhibited rational functions with Julia sets homeomorphic to the Sierpinski carpet.
It is well-known that the residual Julia set is then a connected dense $G_\delta$ subset.

Interestingly, the known examples of Julia sets which have non-empty residual Julia sets are precisely the examples for which no Fatou component has a finite grand orbit.  
Peter M.~Makienko made the following conjecture.

\begin{mconj}
    Let $R:\rsphere \rightarrow \rsphere$ be a rational function.
    The Julia set $J(R)$ has buried points if and only if there is no completely invariant component of the Fatou set of $R^2$.
\end{mconj}

This conjecture was formulated as a possible entry in Sullivan's dictionary in 1990 \cite{Eremenko:1990} as a parallel to a theorem of Abikoff.
Note that one direction is easily proved with facts already given.
Specifically, if a Fatou component $F$ is completely invariant under $R^2$, then $\partial F$ is also completely invariant, closed, and consists of more than three points, so $\partial F = J(R)$.
The example $z \mapsto \frac 1 {z^2}$ illustrates why one must examine the Fatou set of $R^2$.
However, it is known \cite[Theorem 9.4.3]{Beardon:1991} that a rational map may have at most two completely invariant Fatou components, in which case the Julia set is a simple closed curve.

Makienko's conjecture has received attention in the past, with results being limited by topological considerations.
Morosawa \cite{Morosawa:1997,Morosawa:2000} has proved the conjecture in the cases that $R$ is hyperbolic or subhyperbolic.
Qiao \cite{Qiao:1997lr} has proved this conjecture under the assumption that $J(R)$ is locally connected, as well as in the case that $J(R)$ is not connected.
This extends Morosawa's results, since hyperbolic and subhyperbolic rational functions with connected Julia sets have the property that $J(R)$ is locally connected.

An interesting special case was proved by Sun and Yang \cite{Sun:2003}.
Let $R$ be a rational map with exactly two critical points and of degree at least three.
Then either
\begin{enumerate}
  \item $R$ satisfies Makienko's conjecture, or
  \item $J(R)$ is a Lakes of Wada continuum, and therefore either indecomposable or the union of two indecomposable continua.
\end{enumerate}
In \cite{Childers:2006lr}, the second alternative is improved to $J(R)$ being an indecomposable continuum.  We reproduce the theorem used to make this improvement as Theorem~\ref{thm:indec_julia} below.

This collection of results indicates that any counterexample to the conjecture must have a complicated Julia set.
The main theorem of this paper says that the Julia set must be extremely complicated.

\begin{thma}[Makienko's Conjecture for Decomposable Julia Sets of Rational Maps]
    If $R$ is a rational function such that $J(R)$ has no buried points and $F(R^2)$ has no completely invariant components, then $J(R)$ is an indecomposable continuum.
\end{thma}

Recall that a continuum is \emph{decomposable} if it can be written as the union of two proper subcontinua; otherwise it is \emph{indecomposable}. 
There are no known examples of Julia sets which are indecomposable continua.  
In fact, whether or not there exists a rational function with an indecomposable continuum as its Julia set is a well-known unsolved problem \cite{Mayer:1993}.

\end{section}

\begin{section}{Main Result}

For this section, let $R$ be a rational function which is a counterexample to Makienko's conjecture.  Then
\begin{enumerate}
  \item $J(R)$ is connected \cite{Qiao:1997lr},
  \item $J(R) \neq \rsphere$,
  \item $J(R)$ has no buried points, and
  \item $F(R^2)$ has no completely invariant component.
\end{enumerate}
Therefore, $J(R)$ is a continuum and we will show that it is indecomposable. 

The argument will consist of two main parts.
In Subsection~\ref{part1}, we make use the dynamics of $R$ to prove topological facts about $J(R)$, ultimately that it is a continuum which is \emph{irreducible about a finite set}.
In Subsection~\ref{part2}, we apply a decomposition theorem of Kuratowski to show that the Julia set must contain an indecomposable subcontinuum with interior, and therefore must be indecomposable by \cite{Childers:2006lr}.

\subsection{Consequences of the Dynamics}\label{part1}  
In this section we prove our Irreducibility Theorem~\ref{thm:finitely_irreducible}, from which it immediately follows that a counterexample to Makienko's conjecture cannot be arcwise connected.

\begin{lem}\label{lem:unshielded}
  There exists a periodic Fatou component $U$ such that $\bd U = J(R)$.
\end{lem}

\begin{proof}
  If $\mathcal F$ is the collection of Fatou components, then $\bigcup_{F \in \mathcal F} \bd F = J(R)$, since $J(R)$ has no buried points.
  Note that $\mathcal F$ is countable, so the Baire Category Theorem implies that some $F_0 \in \mathcal F$ has boundary with interior in $J(R)$.
  Since $R$ is topologically exact, there exists $k \in \mathbb N$ such that $R^k(\bd F_0) = J(R)$, so $R^k(F_0)$  has boundary equal to $J(R)$.
  By Sullivan's No Wandering Domains Theorem, $R^k(F_0)$ is eventually periodic.
  Let $U$ be the iterate which is periodic with period $n$, and we see that $\bd U = J(R)$.
\end{proof}

\begin{lem}\label{lem:preimage}
	There exists a Fatou component $V\neq U$ such that $R^n(V)=U$.
\end{lem}

\begin{proof}
  Let $U=U_1$, $U_2, \ldots U_n$ be the Fatou components in the cycle of $U$.
  Suppose that such a $V$ does not exist.
  Then $U$, and hence each $U_i$, is completely invariant under $R^n$.
  Since the Fatou set of a rational function can have at most two completely invariant components \cite[Theorem 9.4.3]{Beardon:1991}, we see that $n=2$.
  Thus we have arrived at the contradiction that $U$ is completely invariant under $R^2$.
\end{proof}

A continuum $K$ is said to be \emph{irreducible about} $S$ provided that $S\subset K$ and no proper subcontinuum of $K$ contains $S$.

\begin{thm}[Irreducibility Theorem]\label{thm:finitely_irreducible}
    Suppose $R:\rsphere \rightarrow \rsphere$ is a rational function of degree at least 2 such that $J(R)$ has no buried points and $F(R^2)$ has no completely invariant component.
    Then $J(R)$ is irreducible about a finite point set.
\end{thm}

\begin{proof}
    By Lemma~\ref{lem:unshielded}, let $U$ be a periodic Fatou component whose boundary is $J(R)$.
    Without loss of generality, suppose in fact that $R(U) = U$.
    Let $V \neq U$ be a preimage of $U$.
    By way of contradiction, suppose $J(R)$ is not irreducible about any finite point set.
    We will arrive at a contradiction by first finding a continuum $K \subset \rsphere$ such that
    \begin{enumerate}
      \item $K$ does not contain $J(R)$,
      \item $K$ contains all critical values of $J(R)$, and
      \item $K$ does not separate $U$.
    \end{enumerate}
    We will use $K$ to construct arcs $A_1$ and $A_2$.
    The arc $A_1$ will have the property that its endpoints lie in $U$ and $V$, and the arc $A_2$ will be in the Fatou set connecting the endpoints of $A_1$, contradicting the fact that $U$ and $V$ are distinct Fatou components.

    Let $\widetilde K \subset J(R)$ be a minimal subcontinuum containing a finite set consisting of
    \begin{enumerate}
      \item every critical value in $J(R)$, and
      \item an accessible point on $\bd F$, for every Fatou component $F$ (including $U$) which meets the critical value set.
        (Here, accessible means accessible from within $F$.)
    \end{enumerate}
    By the assumption that $J(R)$ is not irreducible about a finite set, we see that $\widetilde K \neq J(R)$.
    For each Fatou component $F$ intersecting the critical value set, there exists an arc in $\overline F$ which contains $C_R \cap F$ and meets $J(R)$ exactly in a point of $\widetilde K \cap \bd F$.
    Let $K$ denote the continuum composed of $\widetilde K$ and these arcs.
    Notice that $K \cap J(R) = \widetilde K \subsetneq J(R)$, and that $K$ does not separate $U$, as required.

    Let $Q'$ be an evenly covered open subset of $\rsphere \setminus K$ which intersects $J(R)$.
    Since $R(\bd V) = J(R)$, there is a component $Q$ of $R^{-1}(Q')$ which intersects $\bd V$.  
    Note that $Q$ intersects $U$ as well, since $\bd U = J(R)$.
    Let $A_1 \subset Q$ be an arc joining a point of $Q \cap U$ to a point of $Q \cap V$.
    Since $U$ and $V$ each map to $U$, and $R|_{Q}$ is a homeomorphism, $R(A_1) \subset Q'$ is an arc whose endpoints lie in $U$.
    Connect the endpoints of $R(A_1)$ with an arc $A_2' \subset U \setminus K$.

    Note that, since $C' = R(A_1) \cup A_2'$ is disjoint from $K$, $C'$ is contained in a component of $\rsphere \setminus K$.
    That component is a simply connected neighborhood of $C'$ containing no critical values, so each component of $R^{-1}(C')$ maps homeomorphically onto $C'$.
    Let $C$ denote the component of $R^{-1}(C')$ which contains $A_1$.

    Since $R|_{C}$ is a homeomorphism, we have that $C$ is the union of $A_1$ and a preimage $A_2$ of $A_2'$.
    Since $A_2' \subset U$ and $F(R)$ is completely invariant, $A_2 \subset F(R)$.
    However, $A_2$ joins the endpoints of $A_1$, which lie in $U$ and $V$.
    This contradicts that $U$ and $V$ are distinct Fatou components.
\end{proof}

Theorem~\ref{thm:finitely_irreducible} is a strong topological constraint on the nature of $J(R)$.
For instance, any locally connected or path connected continuum which is irreducible about a finite set is necessarily a tree.
This recovers the results of Qiao \cite{Qiao:1997lr} and Morosawa \cite{Morosawa:1997, Morosawa:2000} that $J(R)$ cannot be locally connected.
Furthermore, this additionally proves that $J(R)$ cannot be arcwise connected.

\subsection{Applying Kuratowski's Decomposition}\label{part2}
The class of continua which are irreducible about a finite set contains the class of indecomposable continua. 
However, being an indecomposable continuum is generally much stronger than being irreducible about a finite set. We will investigate the relationship in the present context.

We first state a useful theorem from \cite{Childers:2006fk}.
See also \cite{Rogers:1998}, where a similar theorem first appeared, and \cite{Childers:2006lr} for variants. 

\begin{thm}[ \cite{Childers:2006fk}]\label{thm:nowhere_dense_forever}
    Let $R : \rsphere \rightarrow \rsphere$ be a rational function.
    Let $Y$ be a compact subset of $\rsphere$ and let $X$ be a nowhere dense compact subset of $Y$.
    Then $R(X)$ is nowhere dense in $R(Y)$.
\end{thm}

A crucial topological characteristic of $J(R)$ follows immediately.

\begin{cor}\label{cor:nowhere_dense_forever}
  $\Int_{J(R)}(\bd V) \neq \emptyset$.
\end{cor}

In the present situation, we have more topological information about the continuum $J(R)$ -- it is the boundary of one complementary domain $U$, and another complementary domain $V$ has non-empty interior in its boundary, relative to $J(R)$.
These facts are sufficient to prove that $J(R)$ must contain an indecomposable subcontinuum with non-empty interior, and therefore is itself an indecomposable continuum.

\begin{defn}
    Suppose $X \subset \rsphere$ is the boundary of two connected open sets $U, V \subset \rsphere$.
    Then $X$ is \emph{monostratic} if it is the countable union of indecomposable subcontinua and nowhere dense subcontinua.
\end{defn}

Kuratowski \cite{Kuratowski:1928fk} proved that a continuum which is the common boundary of two regions admits a monotone decomposition to a circle if and only if it is not monostratic.
Each element of the decomposition is the countable union of indecomposable subcontinua and nowhere dense subcontinua of the common boundary.
Thus, any element of the decomposition which contains interior necessarily contains an indecomposable continuum.

\begin{lem}
    Let $U$ be a simply connected open subset of $\rsphere$ with nondegenerate boundary.
    Let $V \neq U$ be a complementary domain of $\bd U$.  Then $\partial V$ is the common boundary of at least two regions:  $V$ and $\operatorname{Comp}(\rsphere \setminus \overline V, U)$, the component of $\rsphere \setminus \overline V$ which contains $U$.
\end{lem}

\begin{proof}
    By the Boundary Bumping Theorem,
    \[\bd\operatorname{Comp}(\rsphere \setminus \overline V, U) \subset \bd V.\]
    Also, observe that
    \[\bd V \subset \overline U \subset \overline{\operatorname{Comp}(\rsphere \setminus \overline V, U)},\] and
    \[\bd V \subset \overline V \subset \overline{\rsphere \setminus \operatorname{Comp}(\rsphere \setminus \overline V, U)}.\]
    These three inclusions imply that $\bd \operatorname{Comp}(\rsphere \setminus \overline V, U) = \bd V$.
\end{proof}

\begin{lem}\label{lem:monostratic}
    Let $U \subset \rsphere$ be a simply connected open set with nondegenerate boundary.
    Let $V \neq U$ be a complementary domain of $\bd U$ such that $\Int_{\bd U}(\bd V) \neq \emptyset$.
    Then, if $\bd V$ is a monostratic continuum, it contains an indecomposable subcontinuum with interior in $\bd U$
\end{lem}

\begin{proof}
    If $\partial V$ is monostratic then by definition $\bd V = \bigcup_{i=1}^\infty X_i \cup \bigcup_{i=1}^\infty A_i$, where each $A_i$ is a nowhere dense subcontinuum of $\bd V$ and each $X_i$ is an indecomposable continuum.
    Note that $\Int_{\bd U}(\bd V) = \bigcup_{i=1}^\infty (X_i \cap \Int_{\bd U}(\bd V)) \cup \bigcup_{i=1}^\infty(A_i \cap \Int_{\bd U}(\bd V))$ is a Baire space, so it is not the union of a countable collection of nowhere dense closed subsets of $\bd U$.
    Since each $A_i$ is nowhere dense in $\bd V$,  and thus nowhere dense in $\bd U$, some $X_i$ has interior in $\bd U$.
\end{proof}

\begin{lem}\label{components}
    If $X$ is a continuum irreducible about $\{a_1, \ldots, a_N\}$, then for each subcontinuum $Y \subset X$ the set $X \setminus Y$ has at most $N$ components.
\end{lem}

\begin{proof}
    If $Y \subset X$ is a continuum and $X \setminus Y$ has more than $N$ components, then at most $N$ components of $X \setminus Y$ can contain points of $\{a_1, \ldots, a_N\}$, since components are disjoint.
   By the Boundary Bumping Theorem, the union of $Y$ with those components is then a proper subcontinuum of $X$ containing $\{a_1, \ldots, a_N\}$, contradicting that $X$ is irreducible about $\{a_1, \ldots, a_N\}$.
\end{proof}

\begin{lem}\label{thm:indec_subcontinuum}
    Let $X = \bd U$ for some simply connected open subset of $\rsphere$ with nondegenerate boundary.
    Let $V \neq U$ be a complementary domain of $X$ such that $\Int_X(\bd V) \neq \emptyset$.
    If $X$ is irreducible about a finite set and $\bd V$ is not monostratic, then $\bd V$ contains an indecomposable subcontinuum with interior in $X$.
\end{lem}

\begin{proof}
    Suppose $X$ is irreducible about a finite set $A$ of cardinality $N$.
    By Lemma~\ref{components}, no subcontinuum of $X$ can separate $X$ into more than $N$ components.
    In particular, $X \setminus \bd V$ consists of connected components $X_1, \ldots, X_k$, where $k \le N$.
    Choose $x_i \in \overline{X_i} \cap \bd V$, and let $B$ denote the finite set
    \[ B = \{x_1, \ldots, x_k\} \cup (A \cap \partial V).\]
    Let $m:\partial V \rightarrow \ucirc$ be the monotone map guaranteed by Kuratowski.
    We now consider two cases: Does $m^{-1}(m(B))$ have interior in $X$?

    Suppose first that it does.
    Then $m^{-1}(m(B))$ is the finite union of the continua $\{m^{-1}(m(x)) \, \mid \, x \in B\}$ which contains an open subset of $X$.
    The Baire Category Theorem implies that there exists $p \in B$ such that $M=m^{-1}(m(p))$ has interior in $X$.
    In particular, $M$ has interior in $\bd V$, and according to the Kuratowski decomposition, $M=(\bigcup_{i=1}^\infty M_i)\cup(\bigcup_{i=1}^\infty C_i)$ where each $M_i$ is indecomposable, and each $C_i$ is a continuum of condensation of $X$.
    An argument analagous to that in the proof of Lemma~\ref{lem:monostratic} shows that some $M_i$ has interior in $X$, and the theorem is proved.

    Otherwise, $m^{-1}(m(B))$ is nowhere dense in $X$. 
    Set 
    \[K = \bd V \setminus m^{-1}(m(B)),\]
    which has non-empty interior in $X$ since $\bd V$ does.
    The set $\ucirc \setminus m(B)$ consists of finitely many disjoint open intervals $I_1, \ldots, I_n$, and $K = m^{-1}(I_1) \cup \ldots \cup m^{-1}(I_n)$.
    Since $K$ contains an open subset of $X$, some $m^{-1}(I_q)$ contains an open subset of $X$.
    Notice that $I_q^c = \ucirc \setminus I_q$ is a closed interval containing $m(B)$. Since $m$ is monotone,  $m^{-1}(I_q^c)$ is a subcontinuum  of $\bd V$ containing $B$.

    Set
        \[ X' = \overline{X_1 \cup \ldots \cup X_k} \cup m^{-1}(I_q^c).\]
    Notice that $X'$ is compact since $I_q^c$ is closed.
    That $m(B) \subset I_q^c$ implies two important facts about $X'$.
    First, since $\{x_1, \ldots, x_k\} \subset m^{-1}(I_q^c)$, we see that $X'$ is connected, since each $\overline{X_i}$ is a connected set which meets the connected set $m^{-1}(I_q^c)$.
    Second, $X' \neq X$, since $X_1 \cup \ldots \cup X_k \cup m^{-1}(I_q^c)$ does not contain $m^{-1}(I_q)$ (as $m^{-1}(I_q)$ has interior in $X$).
    Thus, $X'$ is a proper subcontinuum of $X$ containing $B$, contradicting the fact that $X$ is irreducible about that set.
\end{proof}

\begin{thm}\label{lem:irreducible_implies_indecomposable}
    Let $X = \bd U$ for a simply connected open subset $U$ of $\rsphere$ with nondegenerate boundary.
    Let $V \neq U$ be a complementary domain of $X$ such that $\Int_X(\bd V) \neq \emptyset$.
    If $X$ is irreducible about a finite set, then $\bd V$ contains an indecomposable continuum with interior in \nolinebreak $X$.
\end{thm}

\begin{proof}
    If $\bd V$ is monostratic, then this follows from Lemma~\ref{lem:monostratic}.
    If $\bd V$ is not monostratic, then this follows from Lemma~\ref{thm:indec_subcontinuum}.
\end{proof}

\begin{thm}[Theorem 2.4 from \cite{Childers:2006lr}]\label{thm:indec_julia}
    Suppose $J \neq \rsphere$ is the connected Julia set of a rational map of degree at least two and that $J$ has no buried points.
    If $K$ is an indecomposable subcontinuum of $J$ with nonempty interior in $J$, then $K = J$.
\end{thm}

\begin{cor}[Makienko's Conjecture for Decomposable Julia Sets of Rational Maps]
    If $R$ is a rational function such that $J(R)$ has no buried points and $F(R^2)$ has no completely invariant components, then $J(R)$ is an indecomposable continuum.
\end{cor}

\begin{proof}
    By Lemma~\ref{lem:unshielded}, there is a Fatou component $U$ so that $\bd U = J(R)$.
    By Lemma~\ref{lem:preimage}, let $V \neq U$ be a preimage of $U$.
    Corollary~\ref{cor:nowhere_dense_forever} gives that $\bd V$ has interior in $J(R)$.
    Theorem~\ref{thm:finitely_irreducible} implies that $J(R)$ is irreducible about some finite set.
    Theorem~\ref{lem:irreducible_implies_indecomposable} shows that $J(R)$ contains an indecomposable subcontinuum with interior, and Theorem~\ref{thm:indec_julia} shows that $J(R)$ is indecomposable.
\end{proof}
\end{section}

\begin{section}{Conclusions}
It has now been demonstrated that any counterexample to Makienko conjecture must be exceedingly complicated.  In fact, it is not known if the Julia set of a rational function can be as complicated as required by the conclusion of our theorem.  We restate here the question which appears in \cite{Childers:2006fk, Childers:2006lr, Mayer:1993, Rogers:1998, Sun:2003}.
\begin{ques}
  Can the Julia set of a rational function be an indecomposable continuum?  Can the Julia set of a rational function contain an indecomposable subcontinuum with interior?
\end{ques}
Note that Theorem~\ref{thm:indec_julia} indicates that the answers to these two questions are the same for polynomial Julia sets and for rational functions whose Julia set contains no buried points.

Several authors have extended the study of Makienko's conjecture to transcendental and meromorphic functions.
Dom\'inguez and Fagella \cite{Dominguez:ta} survey the current state of affairs in this direction.
Also see Ng, Zheng, and Choi \cite{Ng:2006lr}, who prove Makienko's conjecture for locally connected Julia sets of certain meromorphic functions.
We do not attempt to apply these techinques to the meromorphic case.
In particular, our approach makes use of Theorem~\ref{thm:nowhere_dense_forever}, which fails dramatically for entire and meromorphic functions.
A particular point of difficulty is summarized in the following question.
\begin{ques}
  Let $f$ be a transcendental meromorphic function.  If $V$ is a Fatou component such that $\bd V$ is nowhere dense in $J(f)$, does it follow that $\bd f(V)$ is also nowhere dense in $J(f)$?
\end{ques}
It is also worth noting that our proof makes use of the non-existence of wandering domains for rational functions, which in general is not true for transcendental functions.
\end{section}
\bibliographystyle{annotate}

\end{document}